\documentclass[a4paper,11pt]{article}
\usepackage{indentfirst,latexsym,bm,color}
\usepackage{amsmath,amssymb,amsfonts}

\usepackage[top=1in,left=1in,right=1in]{geometry}

\makeatletter

\@addtoreset{equation}{section} \makeatother

\begin{document}
\newtheorem{dingli}{Theorem}[section]
\newtheorem{dingyi}[dingli]{Definition}
\newtheorem{tuilun}[dingli]{Corollary}
\newtheorem{zhuyi}[dingli]{Remark}
\newtheorem{yinli}[dingli]{Lemma}
\title{ The absolute values and support projections for a class of operator matrices
involving idempotents
\thanks{ This work was supported by NSF of
China (Nos: 11671242, 11571211) and the
Fundamental Research Funds for the Central Universities
(GK201801011).}}
\author{  Yuan Li\thanks{E-mail address:
 liyuan0401@aliyun.com}, \ Xiaomei Cai, \ Shuaijie Wang }
\date{} \maketitle\begin{center}
\begin{minipage}{16cm}
{ \small $$    \ \ \   \ College \ of \ Mathematics \ and \
Information \ Science,\ Shaanxi \ Normal \ University, $$ $$ Xi'an,
710062, People's \ Republic \ of\ China. $$ }
\end{minipage}
\end{center}
 \vspace{0.05cm}
\begin{center}
\begin{minipage}{16cm}
\ {\small {\bf Abstract }
 Let $\lambda\in \mathbb{R},$ $\mu\in \mathbb{R}$ and $B$ be a linear bounded operator from a Hilbert space $\mathcal{K}$ into another Hilbert space $\mathcal{H}.$ In this paper, we consider the formulas of the absolute value $|Q_{\lambda,\mu}|,$ where $Q_{\lambda,\mu}$ with respect to the decomposition  $\mathcal{H}\oplus\mathcal{K}$ have the operator matrix
  form $Q_{\lambda,\mu}:=\left(\begin{array}{cc}\lambda I&B\\B^*&\mu I\end{array}\right).$ Then the positive part and the support projection of  $Q_{\lambda,0}$  are obtained. Also, we characterize the symmetry $J$ such that a projection $E$ is the $J$-projection.
In particular, the minimal element of the set of all symmetries $J$ with the property $JE\geqslant0$ is described.
\\}
\endabstract
\end{minipage}\vspace{0.10cm}
\begin{minipage}{16cm}
\end{minipage}
\end{center}
\begin{center} \vspace{0.01cm}
\end{center}

\section{ Introduction}
Let $\mathcal{H}$ and $\mathcal{K}$ be complex Hilbert spaces,
$\mathcal{B(K,H)}$ be the set of all bounded linear operators from $\mathcal{K}$ into $\mathcal{H}.$
For an operator $A\in \mathcal{B(H,K)},$ we denote by $|A|$ the absolute value of operator $A,$
that is $|A|=(A^{\ast} A)^{\frac{1}{2}},$ where $A^*$ is the adjoint operator of $A.$
We write $A\geqslant0$ if $A$ is a positive operator, meaning $\langle Ax,x\rangle \geqslant 0,$
where we denote by $\langle,\rangle$ the inner product of $\mathcal{H}.$  Also, denote by $\mathcal{B({H})}^{+}$ the set of all positive bounded linear operators on $\mathcal{H}.$ If $A \in\mathcal{B({H})}^{+},$
then $A^{\frac{1}{2}} $ denotes the positive square root of $A.$
For an operator $T\in \mathcal{B(H,K)},$ we use $N(T),$ $R(T)$ and $\overline{R(T)}$ to
 denote the null
 space, the range of $T,$ and the closure of $R(T),$  respectively.
 If $\mathcal{M}$ is a closed subspace of $\mathcal{H},$ then
 $P_{\mathcal{M}}$ denotes the orthogonal projection onto $\mathcal{M}.$
Particularly, we use $P_{A}$ to denote the orthogonal projection onto
$\overline{R(A)}.$ Also, $P_{A}$ is said to be the support projection of $A$ and
 $E_A(\cdot)$ denotes the spectral measure of $A$
 for a self-adjoint operator $A\in\mathcal{B(H)}.$

As usually, the operator order (Loewner partial order) relation $A\geqslant B$ between two bounded self-adjoint operators is defined as $A-
B\geqslant 0.$
Let $A^+:=\frac{|A|+A}{2}$ and $A^-:=\frac{|A|-A}{2}$ be the positive and negative parts of a
self-adjoint operator $A\in\mathcal{B(H)}.$ Then the self-adjoint operator
$A$ can be written as $A=A^+-A^-$ with $A^+A^-=0.$

It is well known that every operator $T\in \mathcal{B(H)}$ has a (unique) polar decomposition
 $T = U(T^*T)^\frac{1}{2}$
where $U$ is a partial isometry from $N(T)^\perp$ onto $\overline{R(T)}$
with kernel space $N(T)$(See [5]).
 Also, we write $A\simeq B$ to mean that the operator $A$ and $B$ are unitarily
equivalent $(A=UBU^*)$ and $\mathcal{M}\simeq \mathcal{N}$ denotes
$\mathcal{M}=U\mathcal{N}$ for some unitary operator $U,$ where
$\mathcal{M}$ and $\mathcal{N}$ are closed subspaces of $\mathcal{H}.$

An operator $J\in \mathcal{B(H)}$ is said to be a symmetry if $J=J^*=J^{-1}.$
In this case, the operators $J^+:=\frac{I+J}{2}$ and $J^-:=\frac{I-J}{2}$ are
mutually annihilating orthogonal projections. Thus there is a bijection between
the set of all orthogonal projections and all symmetries. If $J$ is a non-scalar symmetry, then
 the form \[[x,y]:=\langle Jx, y\rangle \qquad(x,y\in \mathcal{H})\] is called an indefinite metric and
 $\mathcal{H}$ with $[.,.]$ is said to be a Krein space.

Let us denote by $\mathcal{B(H)}^{Id}$ the set of all projections (idempotents) on $\mathcal{B(H)}.$
 If $E\in \mathcal{B(H)}^{Id},$ then $E$ can be written as a $2\times2$ operator matrix: \begin{equation} E=\left(\begin{array}{cc}I&E_1\\0&0\end{array}\right): R(E)\oplus R(E)^{\perp}\rightarrow R(E)\oplus R(E)^{\perp}\end{equation}
where $E_1\in {\mathcal{B}}(R(E)^{\perp}, R(E)).$ Thus  \begin{equation} E+E^*=\left(\begin{array}{cc}2I&E_1\\E_1^*&0\end{array}\right).\end{equation}
Also, a projection $E$ is said to be $J$-projection if $E=JE^{\ast}J.$
Recently, some interesting decomposition properties of projections and $J$-projections were studied in [3,6,10,12]. An exposition of the properties of the Krein spaces can be found in the book by T. Ya. Azizov and I.S. Iokhvidov [1]. In particular, the
problem of the existence of $J$-projections and its properties are studied in [7-9,11].

In general, with respect to the decomposition $ \mathcal{H}\oplus\mathcal{K},$
 a self-adjoint operator $Q\in \mathcal{B(H\oplus K)}$ has the following matrix representation
$$Q=\left(\begin{array}{cc}Q_1&Q_2\\Q_2^*&Q_3\end{array}\right),$$ where $Q_1$ and $Q_3$ are self-adjoint operators and $Q_2\in \mathcal{B(K,H)}.$
However, the expression of $|Q|$ (in terms of some formulas involving $Q_1,$ $Q_2$ and $Q_3$) seems difficult.
In section 2, we mianly consider how to give the expression of $|Q|,$ on the assumpton of
$Q_1=\lambda I$ and $Q_3=\mu I,$ where $\lambda\in\mathbb{R}$ and $\mu\in\mathbb{R}$ are real numbers.
That is, we study the concrete form of $|Q_{\lambda,\mu}|,$
according to the formulas of
$\lambda,\mu,$ and $B,$ where $$Q_{\lambda,\mu}:=\left(\begin{array}{cc}\lambda I&B\\B^*&\mu I\end{array}\right):\mathcal{H}\oplus\mathcal{K}\rightarrow\mathcal{H}\oplus\mathcal{K}.$$
 The motivation of our research $|Q_{\lambda,\mu}|$  is based on calculating $|E+E^*|$ in (1.2).

In section 3, we mainly give the formulas of the absolute values and support projections
for operator matrices, which have analogous forms as (1.2). In section 4,
 we mainly characterize the symmetry $J$ such that a projection $E$ is the $J$-projection.
In particular, the minimal element of the set of all symmetries $J$ with $JE\geqslant0$ is given. That is
$$ min\{J: JE\geqslant0,\ J=J^{\ast}=J^{-1}\}=2P_{(E+E^{\ast})^{+}}-I.$$
Furthermore, using some new methods and techniques,  we recover some theorems
and corollaries of [2, 11].

\section{ Absolute values for a class of operator matrices }

The following lemma is a direct calculation from diagonalization of a self-adjoint matrix.

 {\bf Lemma 1.}  Let $b>0$ and $\mu\in\mathbb{R}.$ Then

$$\left(\begin{array}{cc}
 1+b&(1+\mu)\sqrt{b}\\\sqrt{b}(1+\mu)&\mu^2+b\end{array}\right)^{\frac{1}{2}}=\begin{cases}\left(\begin{array}{cc}
 1&\sqrt{b}\\\sqrt{b}&\mu\end{array}\right),\ \ \ \ \ \ \ \ \ \ \ \ \ \ \ \ \ \ \ \ \qquad \hbox{ if }  \mu\geqslant b & \\
\dfrac{1}{t}\left(\begin{array}{cc}
 2b-\mu+1&(1+\mu)\sqrt{b}\\ \sqrt{b}(1+\mu)&\mu^2-\mu+2b\end{array}\right), \ \ \ \hbox{ if }\mu< b&
\end{cases}$$ where $t=(\mu^2-2\mu+4b+1)^{\frac{1}{2}}.$

We get the following results inspired by Lemma 1 above.

 {\bf Lemma 2.}  Let $A\in \mathcal{B(H)}^{+}, \mu\in \mathbb{R}.$
 Suppose that $M\in\mathcal{B(H\oplus H)}$ has the operator matrix form
$M:=\left(\begin{array}{cc}
I+A&(1+\mu) A^\frac{1}{2}\\(1+\mu) A^\frac{1}{2}&\mu^2I+A\end{array}\right):\mathcal{H\oplus H}\rightarrow \mathcal{H\oplus H}.$

(i) If $\mu\geqslant\|A\|,$ then
$$ M^{\frac{1}{2}} = \left(\begin{array}{cc}
I& A^{\frac{1}{2}}\\ A^{\frac{1}{2}}&\mu I\end{array}\right):\mathcal{H\oplus H}\rightarrow \mathcal{H\oplus H}.$$

(ii) If $ \mu\leqslant0,$ then
$$ M^{\frac{1}{2}}= \left(\begin{array}{cc}
T^{-1}(2A-\mu I+I)&(1+\mu)T^{-1}A^\frac{1}{2}\\(1+\mu)T^{-1}A^\frac{1}{2}&T^{-1}(\mu^2I-\mu I+2A)\end{array}\right):\mathcal{H\oplus H}\rightarrow \mathcal{H\oplus H},$$
where $T=[(\mu-1)^2I+4A]^{\frac{1}{2}}.$

(iii) If $0<\mu<\|A\|$ and $\mathcal{H}_1:=E_A[0,\mu],$ $\mathcal{H}_2:=E_A(\mu,\|A\|],$ then
with respect to the decomposition $ \mathcal{H}\oplus\mathcal{H}=\mathcal{H}_1\oplus\mathcal{H}_2\oplus\mathcal{H}_1\oplus\mathcal{H}_2\simeq
(\mathcal{H}_1\oplus\mathcal{H}_1)\oplus(\mathcal{H}_2\oplus\mathcal{H}_2),$
$$ \begin{array}{rl}M^{\frac{1}{2}}&=\left(\begin{array}{cccc}
I_1&0&A_1^\frac{1}{2}&0\\0&T_2^{-1}(2A_2-\mu I_{2}+I_2)&0&(1+\mu)T_2^{-1}A_2^\frac{1}{2}\\ A_1^\frac{1}{2}&0&\mu I_1&0\\0&(1+\mu)T_2^{-1} A_2 ^\frac{1}{2}&0&T_2^{-1}(\mu^2I_2-\mu I_2+2A_2)
\end{array}\right)\\&\simeq\left(\begin{array}{cc}
I_1&A_1^\frac{1}{2}\\A_1^\frac{1}{2}&\mu I_1\end{array}\right) \oplus  \left(\begin{array}{cc}
T_2^{-1}(2A_2-\mu I_2+I_2)&(1+\mu)T_2^{-1}A_2^\frac{1}{2}\\(1+\mu)T_2^{-1}A_2^\frac{1}{2}&T_2^{-1}(\mu^2I_2-\mu I_2+2A_2)\end{array}\right),\end{array}$$
where $A_1:=A\mid_{\mathcal{H}_1},$ $A_2:=A\mid_{\mathcal{H}_2}$ and $T_2:=[(\mu-1)^2I_2+4A_2]^\frac{1}{2},$
where $I_i$ are identity operators on the subspace of $H_i$ for $i=1,2.$

{\bf Proof.} (i) and (ii) are easy to verify.

(iii) If $0<\mu<\|A\|$ and $\mathcal{H}_1=E_A[0,\mu],$ $\mathcal{H}_2=E_A(\mu,\|A\|],$
then $A_1=A\mid_{\mathcal{H}_1}\leqslant\mu I_1$ and  $A_2=A\mid_{\mathcal{H}_2}\geqslant\mu I_2,$ so
$\|A_1\|\leqslant\mu$ and $T_2=[(\mu-1)^2I_2+4A_2]^\frac{1}{2}$ is invertible.
It is easy to see that
$$\begin{array}{rl}&(2A_2-\mu I_2+I_2)(\mu^2I_2-\mu I_2+2A_2)-(1+\mu)^2A_2\\=&4A_2^2-6\mu A_2+\mu^2A_2+A_2-\mu(\mu-1)^2I_2\\ \geqslant &(\mu-1)^2(A_2-\mu I_2)\\ \geqslant&0.\end{array}$$
Then the positivity of the Schur complement of the following operator matrix implies
$$\left(\begin{array}{cc}
T_2^{-1}(2A_2-\mu I_2+I_2)&(1+\mu)T_2^{-1}A_2^\frac{1}{2}\\(1+\mu)T_2^{-1} A_2 ^\frac{1}{2}&T_2^{-1}(\mu^2I_2-\mu I_2+2A_2)
\end{array}\right)\geqslant0.$$

Obviously, with respect to the decomposition $H=\mathcal{H}_1\oplus\mathcal{H}_2\oplus\mathcal{H}_1\oplus\mathcal{H}_2,$
the block matrix representation of $M$ is
\begin{equation}M=\left(\begin{array}{cccc}
I_1+A_1&0&(1+\mu)A_1^\frac{1}{2}&0\\0&I_2+A_2&0&(1+\mu)A_2^\frac{1}{2}\\ (1+\mu)A_1^\frac{1}{2}&0&\mu^2I_1+A_1&0\\0&(1+\mu) A_2^\frac{1}{2} &0&\mu^2I_2+A_2
\end{array}\right).\end{equation}
 Then Lemma 1 and a direct calculation imply
$$M^{\frac{1}{2}}=\left(\begin{array}{cccc}
I_1&0&A_1^\frac{1}{2}&0\\0&T_2^{-1}(2A_2-\mu I_2+I_2)&0&(1+\mu)T_2^{-1}A_2^\frac{1}{2}\\ A_1^\frac{1}{2}&0&\mu I_1&0\\0&(1+\mu)T_2^{-1} A_2 ^\frac{1}{2}&0&T_2^{-1}(\mu^2I_2-\mu I_2+2A_2)
\end{array}\right),$$
so with respect to the decomposition $H=(\mathcal{H}_1\oplus\mathcal{H}_1)\oplus(\mathcal{H}_2\oplus\mathcal{H}_2),$ $$M^{\frac{1}{2}}=\left(\begin{array}{cc}
I_1&A_1^\frac{1}{2}\\A_1^\frac{1}{2}&\mu I_1\end{array}\right) \oplus  \left(\begin{array}{cc}
T_2^{-1}(2A_2-\mu I_2+I_2)&(1+\mu)T_2^{-1}A_2^\frac{1}{2}
\\(1+\mu)T_2^{-1}A_2^\frac{1}{2}&T_2^{-1}(\mu^2I_2-\mu I_2+2A_2)\end{array}\right).$$
\qquad $\square$

 Let $D\in \mathcal{B(K,H)}$ with $\dim N(D)=\dim N(D^*)=0.$
 Then the polar decomposition theorem implies that $D^*=V(DD^*)^{\frac{1}{2}},$ where
$V$ is a unitary operator.
Denote $\mathcal{H}_1:=E_{DD^*}[0,\mu]$ and $\mathcal{H}_2:=E_{DD^{\ast}}(\mu,\|D\|^{2}].$
If $\mathcal{K}_1:=V\mathcal{H}_1,$ and $\mathcal{K}_2:=V\mathcal{H}_2,$ then $V$ has the operator matrix form
 \begin{equation} V=\left(\begin{array}{cc}V_{11}&0\\0&V_{22}\end{array}\right):
\mathcal{H}_1\oplus\mathcal{H}_2\longrightarrow\mathcal{K}_1\oplus\mathcal{K}_2. \end{equation}
Also, define $A_1:=DD^{\ast}\mid _{\mathcal{H}_1}$ and $A_2:=DD^{\ast}\mid _{\mathcal{H}_2}.$ Then $D_1^{\ast}:=V_{11}A_1^\frac{1}{2}$ and $D_2^{\ast}:=V_{22} A_2^\frac{1}{2}$ satisfy
  that $D_1^{\ast}\in {\mathcal{B}}({\mathcal{H}_1,\mathcal{K}_1}),$ $D_2^{\ast}\in {\mathcal{B}}({\mathcal{H}_2,\mathcal{K}_2})$ and
 \begin{equation}D^{\ast}=V(DD^{\ast})^{\frac{1}{2}}=\left(\begin{array}{cc}
V_{11}&0\\0&V_{22}\end{array}\right)\left(\begin{array}{cc}
 A_1^\frac{1}{2}&0\\0& A_2^\frac{1}{2}\end{array}\right)=D_1^*\oplus D_2^{\ast}. \end{equation}

 {\bf Lemma 3.} Let $\mu\in \mathbb{R}$ and $D\in \mathcal{B(K,H)}$ with $\dim N(D)=\dim N(D^{\ast})=0.$ If $0<\mu<\|D\|^{2}$ and if $\mathcal{H}_i,$ $\mathcal{K}_i$ and $D_i$ are the same as in (2.2) and (2.3) for $i=1,2,$ then with respect to the space decomposition
$\mathcal{H}\oplus\mathcal{K}=(\oplus_{i=1}^2{\mathcal{H}_i})\bigoplus({\oplus_{i=1}^2} {\mathcal{K}}_i)\simeq\bigoplus_{i=1}^2({\mathcal{H}_i}{\oplus}{\mathcal{K}_i)},$

$$\begin{array}{rl}&\begin{vmatrix}\left(\begin{array}{cc} I & D \\ D^{\ast}&\mu I\end{array}\right)\end{vmatrix}
\\=&\left(\begin{array}{cccc}I_1&0&D_1&0\\0&T_2^{-1}(2D_2D_2^{\ast}-\mu I_2+I_2)&0&(1+\mu)T_2^{-1}D_2\\D_1^{\ast}&0&\mu I_3&0
\\0&(1+\mu)D_2^{\ast}T_2^{-1}&0&(\mu^2 I_4-\mu I_4+2D_2^{\ast}D_2)
S_2^{-1}\end{array}\right)
\\ \simeq&\left(\begin{array}{cc} I_1 & D_1 \\ D^{\ast}_1&\mu I_3\end{array}\right) \oplus \left(\begin{array}{cc} T_2^{-1}(2D_2D_2^{\ast}-\mu I_2+I_2)&(1+\mu)T_2^{-1}D_2\\(1+\mu)D_2^{\ast}T_2^{-1}
&(\mu^2 I_4-\mu I_4+2D_2^{\ast}D_2)
S_2^{-1} \end{array}\right).\end{array}$$
where  $I_i$ (for $i=1,2)$ are identity operators on the subspace of ${\mathcal{H}_i},$
$I_i$ (for $i=3,4)$ are identity operators on the subspace of ${\mathcal{K}}_{i-2},$ $T_2=[(\mu-1)^2I_2+4D_2D_2^{\ast}]^{\frac{1}{2}}$
and $S_2=[(\mu-1)^2I_4+4D_2^{\ast}D_2]^{\frac{1}{2}}.$

{\bf Proof.}
  Define an operator $X$ from $\mathcal{H}\oplus\mathcal{H}=
(\oplus_{i=1}^2{\mathcal{H}_i})\oplus(\oplus_{i=1}^2{\mathcal{H}_i})$ into $\mathcal{H}\oplus\mathcal{K}=(\oplus_{i=1}^2{\mathcal{H}_i})\oplus({\oplus_{i=1}^2} {\mathcal{K}}_i)$ by  $$X:=diag(I_1,I_2,V_{11},V_{22}),$$
where $V_{ii}$ is the same as in (2.2) for $i=1,2.$ Then $X$ is a unitary operator and
$$\left(\begin{array}{cc}I&D\\D^{\ast}&\mu I\end{array}\right)
=X\left(\begin{array}{cccc}I_1&0&A_1^\frac{1}{2}&0\\0&I_2&0& A_2^\frac{1}{2}\\A_1^\frac{1}{2}&0&\mu I_3&0\\0& A_2^\frac{1}{2}&0&\mu I_4\end{array}\right)X^\ast,$$ where
$A_i:=DD^{\ast}\mid _{{\mathcal{H}}_i}$ for $i=1,2,$ which induces
$$\left(\begin{array}{cc}I&D\\D^{\ast}&\mu I\end{array}\right)^{2}=
X\left(\begin{array}{cccc}
I_1+A_1&0&(1+\mu) A_1^\frac{1}{2}&0\\0&I_2+A_2&0&(1+\mu) A_2^\frac{1}{2}\\(1+\mu) A_1^\frac{1}{2}&0&\mu^2I_3+A_1&0\\0&(1+\mu) A_2^\frac{1}{2}&0&\mu^2I_4+A_2\end{array}\right)X^\ast.$$

On the other hand, we have $$T_2^{2}D_2=[(\mu-1)^2I_2+4D_2D_2^{\ast}]D_2=D_2[(\mu-1)^2I_2+4D_2^{\ast}D_2]=D_2S_2^2,$$ which implies $T_2D_2=D_2S_2,$ so $T_2^{-1}D_2=D_2S_2^{-1}.$ Then
$$\begin{array}{rl}&V_{22}T_2^{-1}(\mu^2I_2-\mu I_2+2D_2D_2^{\ast})V_{22}^*
\\=&V_{22}(D_2D_2^*)^\frac{1}{2}(D_2D_2^*)^{-1}T_2^{-1}(\mu^2I_2-\mu I_2+2D_2D_2^{\ast})(D_2D_2^*)^\frac{1}{2}V_{22}^*\\=&D_2^{\ast}(D_2D_2^{\ast})^{-1}(\mu^2I_2-\mu I_2+2D_2D_2^{\ast})T_2^{-1}D_2\\=&
D_2^{-1}(\mu^2I_2-\mu I_2+2 D_2 D_2^{\ast})D_2S_2^{-1}\\=&
D_2^{-1}D_2(\mu^2I_4-\mu I_4+2 D_2^{\ast} D_2)S_2^{-1}\\=&
(\mu^2I_4-\mu I_4+2D_2^{\ast}D_2)S_2^{-1}.\end{array}$$
Thus by Lemma 2 (iii) and a direct calculation, we get
$$\begin{array}{rl}&\begin{vmatrix}\left(\begin{array}{cc}
I&D\\D^{\ast}&\mu I\end{array}\right)\end{vmatrix}=X\left(\begin{array}{cccc}
I_1+A_1&0&(1+\mu) A_1^\frac{1}{2}&0\\0&I_2+A_2&0&(1+\mu) A_2^\frac{1}{2}\\(1+\mu) A_1^\frac{1}{2}&0&\mu^2I_1+A_1&0\\0&(1+\mu) A_2^\frac{1}{2}&0&\mu^2I_2+A_2\end{array}\right)^\frac{1}{2}X^\ast\\=&\left(\begin{array}{cccc}
I_1&0&D_1&0\\0&T_2^{-1}(2D_2D_2^{\ast}-\mu I_2+I_2)&0&(1+\mu)T_2^{-1}D_2\\D_1^{\ast}&0&\mu I_3&0
\\0&(1+\mu)D_2^{\ast}T_2^{-1}&0&V_{22}T_2^{-1}(\mu^2I_2-\mu I_2+2D_2D_2^{\ast})V_{22}^*\end{array}\right)\\=&\left(\begin{array}{cccc}
I_1&0&D_1&0\\0&T_2^{-1}(2D_2D_2^{\ast}-\mu I_2+I_2)&0&(1+\mu)T_2^{-1}D_2\\D_1^{\ast}&0&\mu I_3&0
\\0&(1+\mu)D_2^{\ast}T_2^{-1}&0&(\mu^2I_4-\mu I_4+2D_2^{\ast}D_2)S_2^{-1} \end{array}\right)
 \\ \simeq & \begin{pmatrix} I_1 & D_1 \\ D^{\ast}_1&\mu I_3\end{pmatrix}
\oplus \begin{pmatrix} T_2^{-1}(2D_2D_2^{\ast}-\mu I_2+I_2)&(1+\mu)T_2^{-1}D_2\\
(1+\mu)D_2^{\ast}T_2^{-1}&(\mu^2I_4-\mu I_4+2D_2^{\ast}D_2)S_2^{-1} \end{pmatrix}.\end{array}$$ \ \ \ \ \ \ $\square$

Let $B\in \mathcal{B(K,H)}.$ It is clear that $B$  has the operator matrix form
$$ B=\left(\begin{array}{cc}
\widetilde{B}&0\\0&0\end{array}\right):
N(B)^{\perp}\oplus N(B)\rightarrow\overline{R(B)}\oplus R(B)^{\perp},$$
where $ \widetilde{B}\in B(N(B)^{\perp},\ \overline{R(B)})$ is injective and has dense range.
Then with respect to the space decomposition $\overline{R(B)}\oplus R(B)^{\perp}\oplus N(B)^{\perp}\oplus N(B),$
$$\left(\begin{array}{cc}I&B\\B^{\ast}&\mu I\end{array}\right)=\left(\begin{array}{cccc}
I_1&0& \widetilde{B}&0\\0&I_2&0&0\\ \widetilde{B}^{\ast}&0&\mu I_3&0\\0&0&0&\mu I_4\end{array}\right)
 ,$$
which implies
\begin{equation}\begin{vmatrix}\left(\begin{array}{cc}I&B\\B^{\ast}&\mu I\end{array}\right)\end{vmatrix}
\simeq\begin{vmatrix}\left(\begin{array}{cc}
I_1&\widetilde{B}\\\widetilde{B}^{\ast}&\mu I_3\end{array}\right)\end{vmatrix}
\oplus  I_2 \oplus  |\mu| I_4,\end{equation}
where $I_i$ (for $i=1,2,3,4$) are identity operators on the subspaces of $\overline{R(B)},$ $R(B)^{\perp},$
$N(B)^{\perp}$ and $N(B),$ respectively.

{\bf Remark.} Let $0<\mu<\|B\|^{2}$ and $\widetilde{B}^*=\widetilde{V}(\widetilde{B}\widetilde{B}^*)^\frac{1}{2}$
 be the polar decomposition of $\widetilde{B}^*,$
 where $\widetilde{V}$ is the unique unitary operator from $\overline{R(B)}$ onto $N(B)^{\perp}.$
 In a similar way to (2.2) and (2.3), we denote $\widetilde{\mathcal{H}_1}:=E_{\widetilde{B}\widetilde{B}^*}[0,\mu]$ and $\widetilde{\mathcal{H}_2}:=E_{\widetilde{B}\widetilde{B}^{\ast}}(\mu,\|D\|^{2}].$
 Thus $\overline{R(B)}=\widetilde{\mathcal{H}_1}\oplus\widetilde{\mathcal{H}_2}.$

If $\widetilde{\mathcal{K}_1}:=\widetilde{V}\widetilde{\mathcal{H}_1}$ and $\widetilde{\mathcal{K}_2}:=\widetilde{V}\widetilde{\mathcal{H}_2},$ then $N(B)^{\perp}=\widetilde{\mathcal{K}_1}\oplus\widetilde{\mathcal{K}_2}$ and
$\widetilde{V}$ has the operator matrix form
 \begin{equation} \widetilde{V}=\left(\begin{array}{cc}\widetilde{V}_{11}&0\\0&\widetilde{V}_{22}\end{array}\right):
\widetilde{\mathcal{H}_1}\oplus\widetilde{\mathcal{H}_2}\longrightarrow
\widetilde{\mathcal{K}_1}\oplus\widetilde{\mathcal{K}_2}. \end{equation}
Analogously, define $\widetilde{A_1}:=\widetilde{B}\widetilde{B}^{\ast}\mid _{\widetilde{\mathcal{H}_1}}$
and $\widetilde{A_2}:=\widetilde{B}\widetilde{B}^{\ast}\mid _{\widetilde{\mathcal{H}_2}}.$ Then $\widetilde{
B_1}^{\ast}:=\widetilde{V}_{11}\widetilde{A_1}^\frac{1}{2}$ and $\widetilde{B_2}^{\ast}:=\widetilde{V}_{22}\widetilde{A_2}^\frac{1}{2}$ satisfy
  that $\widetilde{B_1}^{\ast}\in {\mathcal{B}}({\widetilde{\mathcal{H}_1},\widetilde{\mathcal{K}_1}}),$ $\widetilde{B_2}^{\ast}\in {\mathcal{B}}(\widetilde{{\mathcal{H}_2}},\widetilde{\mathcal{K}_2})$ and
 $\widetilde{B}^{\ast}=\widetilde{B_1}^*\oplus \widetilde{B_2}^{\ast}.$
 Also, $I_{1i}$ and $I_{3i}$ (for $i=1,2$) are identity operators on the subspaces of
 $\widetilde{\mathcal{H}_i}$ and $\widetilde{\mathcal{K}_i},$
  respectively.

{\bf Lemma 4.} Let $B\in \mathcal{B(K,H)}$ and $\mu\in \mathbb{R}.$ Then

(i) If $\mu\geqslant \|B\|^{2},$ then
$$\begin{vmatrix}\left(\begin{array}{cc}
I&B\\B^{\ast}&\mu I\end{array}\right)\end{vmatrix}=\left(\begin{array}{cc}
I&B\\B^{\ast}&\mu I\end{array}\right):\mathcal{H}\oplus\mathcal{K}\rightarrow
\mathcal{H}\oplus\mathcal{K}.$$

(ii) If $\mu\leqslant0,$ then with respect to the decomposition $\mathcal{H}\oplus\mathcal{K},$ we have

$$\begin{vmatrix}\left(\begin{array}{cc}
I&B\\B^{\ast}&\mu I\end{array}\right)\end{vmatrix}
=\left(\begin{array}{cc}
 T^{-1}(2BB^{\ast}-\mu I+I)&(1+\mu)T^{-1}B\\(1+\mu)B^{\ast} T^{-1}&S^{-1}[(\mu^2-\mu)I+2B^*B]\end{array}\right),$$
where $T:=[(\mu-1)^2I+4BB^{\ast}]^{\frac{1}{2}}$ and $S:=[(\mu-1)^2I+4B^*B]^{\frac{1}{2}}.$

(iii) If $0<\mu<\|B\|^{2},$ then
with respect to the decomposition $\mathcal{H}\oplus\mathcal{K}\simeq(\widetilde{\mathcal{H}_2}\oplus\widetilde{\mathcal{K}_2})\oplus
(\widetilde{\mathcal{H}_1}\oplus\widetilde{\mathcal{K}_1})
\oplus R(B)^\perp\oplus N(B),$ we have

$$\begin{array}{rl}\begin{vmatrix}\left(\begin{array}{cc}
I&B\\B^{\ast}&\mu I\end{array}\right)\end{vmatrix}
\simeq&\left(\begin{array}{cc}
\widetilde{T_2}^{-1}(2\widetilde{B_2}\widetilde{B_2}^{\ast}-\mu I_{12}+I_{12})&(1+\mu)\widetilde{T_2}^{-1}\widetilde{B_2}
\\(1+\mu)\widetilde{B_2}^{\ast}\widetilde{T_2}^{-1}&(\mu^2I_{32}-\mu I_{32}+2\widetilde{B_2}^{\ast}\widetilde{B_2})
\widetilde {S_2}^{-1}\end{array}\right)
 \\ \oplus&\left(\begin{array}{cc}
I_{11}&\widetilde{B_1}\\ \widetilde{B_1}^{\ast}&\mu I_{31}\end{array}\right)\oplus I_2\oplus|\mu|I_4,\end{array}$$
 where $\widetilde{B_i},$ $I_{1i}$ and $I_{3i}$ (for $i=1,2$) are as above Remark,  $\widetilde {T_2}=[(\mu-1)^2I_{12}+4\widetilde{B_2}\widetilde{B_2}^{\ast}]^{\frac{1}{2}}$ and
 $\widetilde {S_2}=[(\mu-1)^2I_{32}+4\widetilde{B_2}^{\ast}\widetilde{B_2}]^{\frac{1}{2}}$

{\bf Proof.} (i) is clear.

(ii) Obviously, $\mu\leqslant0$ implies that $T=[(\mu-1)^2I+4BB^{\ast}]^{\frac{1}{2}}$ and $S=[(\mu-1)^2I+4B^*B]^{\frac{1}{2}}$ are invertible.
Let $U$ be the unique partial isometry such that
$B^{\ast}=U(BB^{\ast})^{\frac{1}{2}},$ $R(U)=\overline{R(B^*)}$ and $R(U^*)=\overline{R(B)}.$
By a direct calculation, we get that   $$\begin{array}{rl}&\left(\begin{array}{cc}
I&B\\B^{\ast}&\mu I\end{array}\right)^2
=\left(\begin{array}{cc}I+BB^{\ast}&(1+\mu)B\\(1+\mu)B^{\ast}&\mu^2I+B^{\ast}B\end{array}\right)
\\=&\left(\begin{array}{cc}I&0\\0&U\end{array}\right)
\left(\begin{array}{cc}I+BB^{\ast}&(1+\mu)|B^{\ast}|\\(1+\mu)|B^{\ast}|&\mu^2I+BB^{\ast}\end{array}\right)
\left(\begin{array}{cc}I&0\\0&U^{\ast}\end{array}\right)+
\left(\begin{array}{cc}0&0\\0&\mu^2P_{B^*}^\perp\end{array}\right)\\=&\left(\begin{array}{cc}I&0\\0&UP_{B}\end{array}\right)
\left(\begin{array}{cc}I+BB^{\ast}&(1+\mu)|B^{\ast}|\\(1+\mu)|B^{\ast}|&\mu^2I+BB^{\ast}\end{array}\right)
\left(\begin{array}{cc}I&0\\0&P_{B}U^{\ast}\end{array}\right)+
\left(\begin{array}{cc}0&0\\0&\mu^2P_{B^*}^\perp\end{array}\right)
\\=&\left(\begin{array}{cc}I&0\\0&U\end{array}\right)
\left(\begin{array}{cc}I+BB^{\ast}&(1+\mu)|B^{\ast}|\\(1+\mu)|B^{\ast}|&\mu^2P_{B}+BB^{\ast}\end{array}\right)
\left(\begin{array}{cc}I&0\\0&U^{\ast}\end{array}\right)+
\left(\begin{array}{cc}0&0\\0&\mu^2P_{B^*}^\perp\end{array}\right),\end{array}$$
since $UU^*=P_{B^*},$ $U(BB^*)U^*=B^{\ast}B$ and $U P_{B}=U.$  Clearly,
$$\left(\begin{array}{cc}I+BB^{\ast}&(1+\mu)|B^{\ast}|
\\(1+\mu)|B^{\ast}|&\mu^2P_{B}+BB^{\ast}\end{array}\right)^{\frac{1}{2}}=
\left(\begin{array}{cc}I&0\\0&P_{B}\end{array}\right)
\left(\begin{array}{cc}I+BB^{\ast}&(1+\mu)|B^{\ast}|\\(1+\mu)|
B^{\ast}|&\mu^2I+BB^{\ast}\end{array}\right)^{\frac{1}{2}},$$
$P_{B^*}^\perp B^*=0$ and $BP_{B^*}^\perp =0.$
Setting $$\widetilde{U}=\left(\begin{array}{cc}I&0\\0&U\end{array}\right):
\mathcal{H\oplus H}\rightarrow\mathcal{H\oplus K},$$ we conclude from Lemma 2 (ii) that
$$\begin{array}{rl}&\begin{vmatrix}\left(\begin{array}{cc}I&B\\B^{\ast}&\mu I\end{array}\right)\end{vmatrix}
\\=&\widetilde{U}
\left(\begin{array}{cc}I+BB^{\ast}&(1+\mu)|B^{\ast}|\\(1+\mu)|B^{\ast}|&\mu^2P_{B}+BB^{\ast}
\end{array}\right)^{\frac{1}{2}}
\widetilde{U}^*+
\left(\begin{array}{cc}0&0\\0&\mu^2P_{B^*}^\perp\end{array}\right)^\frac{1}{2}
\\=&\left(\begin{array}{cc}I&0\\0&UP_{B}\end{array}\right)
\left(\begin{array}{cc}I+BB^{\ast}&(1+\mu)|B^{\ast}|\\(1+\mu)|B^{\ast}|&\mu^2I+BB^{\ast}
\end{array}\right)^{\frac{1}{2}}
\left(\begin{array}{cc}I&0\\0&U^{\ast}\end{array}\right)+
\left(\begin{array}{cc}0&0\\0&|\mu|P_{B^*}^\perp\end{array}\right)
\\=&\widetilde{U}
\left(\begin{array}{cc}
T^{-1}(2BB^{\ast}-\mu I+I)&(1+\mu)T^{-1}|B^{\ast}|\\(1+\mu)|B^{\ast}|T^{-1}&T^{-1}(\mu^2I-\mu I+2BB^{\ast})\end{array}\right)
\widetilde{U}^*+
diag(0,|\mu|P_{B^*}^\perp)
\\=&\left(\begin{array}{cc}
T^{-1}(2BB^{\ast}-\mu I+I)&(1+\mu)T^{-1}B\\(1+\mu)B^{\ast}T^{-1}&UT^{-1}(\mu^2I-\mu I+2BB^{\ast})U^{\ast}+|\mu|P_{B^*}^\perp\end{array}\right).\end{array}$$

Moreover, $B^{\ast}=U(BB^{\ast})^{\frac{1}{2}},$ implies
$(B^{\ast}B)^{\frac{1}{2}}=U(BB^{\ast})^{\frac{1}{2}}U^*,$  so $$(B^{\ast}B)^{\frac{1}{2}}U=U(BB^{\ast})^{\frac{1}{2}}U^*U=U(BB^{\ast})^{\frac{1}{2}}=B^*,$$ which yields
$$\begin{array}{rcl}U[(\mu-1)^2I+4BB^{\ast}]&=&(\mu-1)^2 U+4B^*(BB^{\ast})^{\frac{1}{2}}
\\&=&(\mu-1)^2 U+4(B^*B)^{\frac{1}{2}}B^{\ast}\\&=&
[(\mu-1)^2 +4B^*B]U.\end{array}$$ Thus $$UT^{-1}=U[(\mu-1)^2I+4BB^{\ast}]^{-{\frac{1}{2}}}=[(\mu-1)^2 +4B^*B]^{-{\frac{1}{2}}}U=S^{-1}U,$$  which induces
$$\begin{array}{rl}&UT^{-1}(\mu^2I-\mu I+2BB^{\ast})U^{\ast}+|\mu|P_{B^*}^\perp\\=&
[(\mu-1)^2 +4B^*B]^{-\frac{1}{2}}[(\mu^2-\mu)P_{B^*}+2B^*B]-\mu P_{B^*}^\perp\\=&
S^{-1}[(\mu^2-\mu)P_{B^*}+2B^*B-\mu(1-\mu)P_{B^*}^\perp]\\=&
S^{-1}[(\mu^2-\mu)I+2B^*B],\end{array}$$ since $P_{B^*}^\perp S^2=P_{B^*}^\perp[(\mu-1)^2 +4B^*B]=(1-\mu)^2P_{B^*}^\perp$ implies $ SP_{B^*}^\perp=(1-\mu)P_{B^*}^\perp.$

 (iii) According to equations (2.4),(2.5) and Lemma 3, we know that $$\begin{array}{rcl}\begin{vmatrix}\left(\begin{array}{cc}I&B\\B^{\ast}&\mu I\end{array}\right)\end{vmatrix}
&\simeq& \begin{vmatrix}\left(\begin{array}{cc}
I_1&\widetilde{B}\\ \widetilde{B}^{\ast}&\mu I_3\end{array}\right)\end{vmatrix}
\oplus  I_2 \oplus |\mu| I_4\\&=&\left(\begin{array}{cc}
\widetilde{T_2}^{-1}(2\widetilde{B_2}\widetilde{B_2}^{\ast}-\mu I_{12}+I_{12})&(1+\mu)\widetilde{T_2}^{-1}\widetilde{B_2}
\\(1+\mu)\widetilde{B_2}^{\ast}
\widetilde{T_2}^{-1}&(\mu^2I_{32}-\mu I_{32}+2\widetilde{B_2}\widetilde{B_2}^{\ast})
\widetilde {S_2}^{-1}\end{array}\right)
 \\& \oplus&\left(\begin{array}{cc}
I_{11}&\widetilde{B_1}\\ \widetilde{B_1}^{\ast}&\mu I_{31}\end{array}\right)\oplus I_2\oplus|\mu|I_4. \end{array}$$
 \qquad $\square$

{\bf Theorem 5.} Let $\lambda\in \mathbb{R}, \mu\in  \mathbb{R}$ and  $B\in \mathcal{B(K,H)}.$

 (i) If $\lambda=\mu=0,$ then

$$\begin{vmatrix}\left(\begin{array}{cc}
\lambda I&B\\B^{\ast}&\mu I\end{array}\right)\end{vmatrix}=\left(\begin{array}{cc}
|B^{\ast}|&0\\0&|B|\end{array}\right):\mathcal{H}\oplus\mathcal{K}\rightarrow
\mathcal{H}\oplus\mathcal{K}.$$

(ii) If $\lambda \neq 0$ and $\mu=0,$ then

 $$\begin{vmatrix}\left(\begin{array}{cc}
\lambda I&B\\B^{\ast}&0\end{array}\right)\end{vmatrix}=
\left(\begin{array}{cc}
\frac{T+\lambda^{2} T^{-1}}{2}&\lambda T^{-1}B
\\ \lambda B^{\ast}T^{-1} &2B^{\ast}T^{-1}B \end{array}\right):\mathcal{H}\oplus\mathcal{K}\rightarrow
\mathcal{H}\oplus\mathcal{K},$$
where $T=(\lambda^{2}I+4BB^{\ast})^{\frac{1}{2}}.$

(iii) If $\mu\lambda\geqslant\|B\|^{2}>0,$ then

$$ \begin{vmatrix}\left(\begin{array}{cc}
\lambda I&B\\B^{\ast}&\mu I\end{array}\right)\end{vmatrix}=\frac{|\lambda|}
{\lambda}\left(\begin{array}{cc}
\lambda I&B \\ B^{\ast}&\mu I\end{array}\right):\mathcal{H}\oplus\mathcal{K}\rightarrow
\mathcal{H}\oplus\mathcal{K}.$$

(iv) If $\mu\lambda<0,$ then with respect to the decomposition $\mathcal{H}\oplus\mathcal{K}:$
$$\begin{vmatrix}\left(\begin{array}{cc}
\lambda I&B\\B^{\ast}&\mu I\end{array}\right)\end{vmatrix}
=\left(\begin{array}{cc}
T^{-1}(2BB^{\ast}-\lambda \mu I+\lambda^{2}I)
&(\lambda+\mu)T^{-1}B
\\ (\lambda+\mu)B^*T^{-1} &S^{-1}(\mu^2I-\lambda \mu I+2B^{\ast}B)\end{array}\right),$$
where $T=[(\mu-\lambda)^2 I+4BB^{\ast}]^\frac{1}{2}$ and $S=[(\mu-\lambda)^2 I+4B^{\ast}B]^\frac{1}{2}.$

(v) If $0<\mu\lambda<\|B\|^{2},$ then
with respect to the decomposition $\mathcal{H}\oplus\mathcal{K}\simeq R(B)^\perp\oplus N(B)\oplus(\widetilde{\mathcal{H}_1}\oplus\widetilde{\mathcal{K}_1})
\oplus(\widetilde{\mathcal{H}_2}\oplus\widetilde{\mathcal{K}_2}),$ we have
$$\begin{array}{rl}\begin{vmatrix}\left(\begin{array}{cc}
\lambda I&B\\B^{\ast}&\mu I\end{array}\right)\end{vmatrix}&\simeq |\lambda| I_2\oplus|\mu|I_4\oplus|\lambda|\left(\begin{array}{cc}
I_{11}&\frac{\widetilde{B_1}}{\lambda}\\ \frac{\widetilde{B_1}^{\ast}}{\lambda}&\frac{\mu}{\lambda}I_{31}\end{array}\right)
\\& \oplus \left(\begin{array}{cc}
\widetilde{T}^{-1}(2\widetilde{B_2} \widetilde{B_2}^{\ast}-{\lambda}\mu I
+{\lambda^{2}}I)&
({\lambda}+\mu)\widetilde{T}^{-1}\widetilde{B_2} \\
({\lambda}+\mu)\widetilde{B_2}^{\ast} \widetilde{T}^{-1}&(\mu^2I_{32}-\lambda\mu I_{32}+2\widetilde{B_2}^{\ast}\widetilde{B_2})
\widetilde {S_2}^{-1}\end{array}\right),\end{array}$$
where $\widetilde{B_i}^{\ast}$ is the same to that of Lemma 4 (for $i=1,2$), $\widetilde{T}=[(\mu-\lambda)^2 I_{12}+4\widetilde{B_2}\widetilde{B_2}^{\ast}]^\frac{1}{2}$ and $\widetilde{S}=[(\mu-\lambda)^2 I_{32}+4\widetilde{B_2}^{\ast}\widetilde{B_2}]^\frac{1}{2}.$

(vi) If $\lambda=0$ and $\mu\neq 0,$ then

 $$\begin{vmatrix}\left(\begin{array}{cc}
0&B\\B^{\ast}&\mu I\end{array}\right)\end{vmatrix}=
\left(\begin{array}{cc}
2BS^{-1}B^{\ast}&{\mu}B S^{-1}
\\ {\mu}S^{-1}B^{\ast} &\frac{S+{\mu}^{2}S^{-1}}{2}\end{array}\right):\mathcal{H}\oplus\mathcal{K}\rightarrow
\mathcal{H}\oplus\mathcal{K},$$
where $S=(\mu^{2}I+4B^*B)^{\frac{1}{2}}.$

{\bf Proof.} (i) and (iii) are direct calculations. It is easy to verify that Lemma 4 (ii) implies (ii).
(iv) follows from Lemma 4 (ii) and the fact that if $\lambda\neq 0,$ then \begin{equation}\begin{vmatrix}\left(\begin{array}{cc}
\lambda I&B\\B^{\ast}&{\mu}I \end{array}\right)\end{vmatrix}=|\lambda|
\begin{vmatrix}\left(\begin{array}{cc}
I&\frac{B}{\lambda }\\ \frac{B^{\ast}}{\lambda}&\frac{{\mu}}{\lambda}I\end{array}\right)\end{vmatrix}.\end{equation}
Using equation (2.6) and Lemma 4 (iii), we get (v).
 (vi) is obvious from (ii) and the relation $\begin{vmatrix}\left(\begin{array}{cc}
0&B\\B^{\ast}&\mu I\end{array}\right)\end{vmatrix}\simeq \begin{vmatrix}\left(\begin{array}{cc}
 \mu I&B^{\ast}\\B&0\end{array}\right)\end{vmatrix}.$ \qquad $\square$

\section{Support projections for a class of operator matrices}

In this section, we will give the expressions of the positive parts and support projections for operator matrices \begin{equation}S_\lambda:=Q_{\lambda,0}=\left(\begin{array}{cc}\lambda I&B\\B^{\ast}&0\end{array}\right):
\mathcal{H}\oplus\mathcal{K}\rightarrow
\mathcal{H}\oplus\mathcal{K},\end{equation}  where $\lambda\in\mathbb{R}$ and $B\in \mathcal{B(K,H)}.$ The following lemma
follows from diagonalization of a self-adjoint matrix ([4]).

{\bf Lemma 6.} Let $a>0.$ Then $C:=\left(\begin{array}{cc}1&2\sqrt{a}
\\2\sqrt{a}&4a\end{array}\right)$ is positive
and the support projection of $C$ is
$$P_{C}= \left(\begin{array}{cc}\dfrac{1}{1+4a}&\dfrac{2\sqrt{a}}{1+4a}
\\\\ \dfrac{2\sqrt{a}}{1+4a}&\dfrac{4a}{1+4a}\end{array}\right).$$

Using Lemma 6,  we get some hints for the following results.

{\bf Lemma 7.}  Let $A\in \mathcal{B(H)}^{+}.$
 If $\widetilde{A}\in B(\mathcal{H}\oplus\mathcal{H})$ has the operator matrix form  $\widetilde{A}=\left(\begin{array}{cc}
 I&2A^{\frac{1}{2}}\\2A^{\frac{1}{2}}&
 4A\end{array}\right),$
then $\widetilde{A}\in B(\mathcal{H}\oplus\mathcal{H})^+$ and the support projection of $\widetilde{A}$ is
$$P_{\widetilde{A}}:=\left(\begin{array}{cc} (I+4A)^{-1}&2A^{\frac{1}{2}}(I+4A)^{-1}\\2A^{\frac{1}{2}}(I+4A)^{-1}&4A(I+4A)^{-1}\end{array}\right).$$

{\bf Proof.} Clearly, $\widetilde{A}\geqslant0.$
It is easy to verify that \begin{equation}P_{\widetilde{A}}^{2}=P_{\widetilde{A}}\geqslant0\ \ \ \ \ \ \ \ \hbox{and} \ \ \ \ \ \ \ \ P_{\widetilde{A}}\widetilde{A}=\widetilde{A}.\end{equation}
Let $x \in\mathcal{H}$ and $y\in\mathcal{H}$ satisfy $\widetilde{A}(x,y)^{t}=0.$ That is
$$\begin{cases}x+2A^{\frac{1}{2}}y
=0&\\2A^{\frac{1}{2}}x+4Ay=0,
\end{cases}$$
then $x=-2A^{\frac{1}{2}}y,$ so
$$P_{\widetilde{A}}(x,y)^{t}=0.$$ Thus
$N(\widetilde{A})\subseteq N(P_{\widetilde{A}})$ yields $\overline{R (\widetilde{A})}\supseteq R(P_{\widetilde{A}}).$
Then by equation (3.2), $\overline{R(\widetilde{A})}= R(P_{\widetilde{A}}),$ which says that
$P_{\widetilde{A}}$ is the support projection of $\widetilde{A}. \ \ \ \ \  \ \ \ \ \ \ \ \square $

{\bf Lemma 8.} Let $A,B\in \mathcal{B(H)}^{+}.$ If $B$ is invertible and $AB=BA,$ then the support projection of $\widetilde{S}:=\left(\begin{array}{cc}
B^2&2A^{\frac{1}{2}}B\\2A^{\frac{1}{2}}B&4A
\end{array}\right):\mathcal{H}\oplus\mathcal{H}
\rightarrow\mathcal{H}\oplus\mathcal{H}$ is $$P_{\widetilde{S}}=\left(\begin{array}{cc}B^2(B^2+4A)^{-1}&2A^{\frac{1}{2}}B(B^2+4A)^{-1}
\\2A^{\frac{1}{2}}B(B^2+4A)^{-1}&4A(B^2+4A)^{-1}
\end{array}\right):\mathcal{H}\oplus\mathcal{H}
\rightarrow\mathcal{H}\oplus\mathcal{H}.$$

{\bf Proof.} Obviously, $AB=BA$ yields $AB^{-1}=B^{-1}A$ and $A^{\frac{1}{2}}B^{-1}=B^{-1}A^{\frac{1}{2}},$ so with respect to the decomposition $\mathcal{H}\oplus\mathcal{H}$
$$S:=\left(\begin{array}{cc}I&2A^{\frac{1}{2}}B^{-1}\\2A^{\frac{1}{2}}B^{-1}&4AB^{-2}\end{array}\right)\geqslant 0.$$
It is easy to see that $\widetilde{S}$ can be written as
$$\widetilde{S}=S\left(\begin{array}{cc}B^2&0\\0&B^2\end{array}\right)
=\left(\begin{array}{cc}B^2&0\\0&B^2\end{array}\right)S,$$ and $diag(B^2,B^2)$ is invertible on the space $\mathcal{H}\oplus \mathcal{H}.$
Then $R(\widetilde{S})=R(S),$ so Lemma 7 implies that with respect to the decomposition $\mathcal{H}\oplus\mathcal{H}$
 $$\begin{array}{rl}P_{\widetilde{S}}=P_{S}&=
\left(\begin{array}{cc}(I+4AB^{-2})^{-1}&2A^{\frac{1}{2}}B^{-1}(I+4AB^{-2})^{-1}
\\2A^{\frac{1}{2}}B^{-1}(I+4AB^{-2})^{-1}&4AB^{-2}(I+4AB^{-2})^{-1}\end{array}\right)
\\\\&=\left(\begin{array}{cc}B^2(B^2+4A)^{-1}&2A^{\frac{1}{2}}B(B^2+4A)^{-1}
\\2A^{\frac{1}{2}}B(B^2+4A)^{-1}&4A(B^2+4A)^{-1}\end{array}\right).\end{array}$$
\qquad $\square$

{\bf Lemma 9.} Let $F\in \mathcal{B(K)}$ and $G\in \mathcal{B(H)}$ be self-adjoint operators.
If there exists a partial isometry operator $U\in \mathcal{B(H,K)}$
such that $F=UGU^{\ast}$ and $U^{\ast}UG=G,$ then $P_{F}=UP_{G}U^{\ast}.$

{\bf Proof.} Denote $Q:=UP_{G}U^{\ast}.$ Then $$Q^2=UP_{G}U^{\ast}UP_{G}U^{\ast}
=UP_{G}U^{\ast}=Q,$$ so $Q$ is an orthogonal projection. Taking $x\in \mathcal{H},$
we get that  $$ Fx=UGU^{\ast}x=UP_{G}U^{\ast}UGU^{\ast}x=QUGU^{\ast}x=QFx,$$ which yields
$R(F)\subseteq R(Q),$ so $ \overline{R(F)}\subseteq R(Q).$

On the other hand, if $Fx=0,$ then $$GU^{\ast}x=U^{\ast}UGU^{\ast}x=U^*Fx=0,$$
which implies that $U^{\ast}x\in N(G)=R(G)^{\perp},$ so $ P_{G}U^{\ast}x=0.$
Hence, $$Qx=UP_{G}U^{\ast}x=0,$$ which induces
$N(F)\subseteq N(Q),$ that is $\overline{R(F)}\supseteq R(Q).$
Thus $\overline{R(F)}=R(Q),$ which says  $P_{F}=Q=UP_{G}U^{\ast}.$  \ \ \ \ \ $ \square$

The following theorem is the main result of this section.

{\bf Theorem 10.} Let $B\in \mathcal{B(K,H)}$ and $S_\lambda$ is the same as (3.1) for $\lambda\in\mathbb{R}.$   Then

(i) If $\lambda\neq0,$ then $$S_{\lambda}^{+}=\frac{1}{2} \left(\begin{array}{cc}
\frac{T_\lambda+\lambda^2T_\lambda^{-1}
+2\lambda I}{2}&B+\lambda T_\lambda^{-1}B\\ B^{\ast}+ \lambda B^*T_\lambda^{-1}&
2B^{\ast}T_{\lambda}^{-1}B\end{array}\right):
\mathcal{H}\oplus\mathcal{K}\rightarrow \mathcal{H}\oplus\mathcal{K},$$ where $T_\lambda=(\lambda^2 I+4BB^{\ast})^{\frac{1}{2}}.$

(ii) $$P_{S_1^+}=\left(\begin{array}{cc}\frac{I+T^{-1}}{2}&T^{-1}B\\B^{\ast}T^{-1}&2B^{\ast}T^{-1}(I+T)^{-1}B\end{array}
\right):
\mathcal{H}\oplus\mathcal{K}\rightarrow \mathcal{H}\oplus\mathcal{K},$$
where $T:=T_1=(I+4BB^{\ast})^{\frac{1}{2}}.$

(iii) If $\lambda>0,$ then $$P_{S_\lambda^+}=\left(\begin{array}{cc}
\frac{I+\lambda T_\lambda^{-1}}{2}
&T_\lambda^{-1}B\\ B^{\ast}T_\lambda^{-1}&
2B^{\ast}T_\lambda^{-1}(\lambda I+T_\lambda)^{-1}B\end{array}
\right):
\mathcal{H}\oplus\mathcal{K}\rightarrow \mathcal{H}\oplus\mathcal{K},$$
where $T_\lambda=(\lambda^2 I+4BB^{\ast})^{\frac{1}{2}}.$

(iv) If $\lambda<0,$ then $$P_{S_\lambda^+}=\left(\begin{array}{cc}
\frac{I+\lambda T_\lambda^{-1}}{2}
&T_\lambda^{-1}B\\ B^{\ast}T_\lambda^{-1}
&\frac{V(I-\lambda T_\lambda^{-1})V^*}{2}\end{array}
\right):
\mathcal{H}\oplus\mathcal{K}\rightarrow \mathcal{H}\oplus\mathcal{K},$$ where $T_\lambda=(\lambda^2 I+4BB^{\ast})^{\frac{1}{2}}$ as above and $V$ is the unique partial isometry such that
$B^{\ast}=V(BB^{\ast})^{\frac{1}{2}}$ with $R(V)=\overline{R(B^*)}$ and $R(V^*)=\overline{R(B)}.$

(v) $$P_{S_0^+}=\frac{1}{2}\left(\begin{array}{cc}
P_B&V^*\\V&P_{B^*}\end{array}
\right):
\mathcal{H}\oplus\mathcal{K}\rightarrow \mathcal{H}\oplus\mathcal{K},$$ where $V$ is the same to (iv).

{\bf Proof.} (i) is immediate from the fact $S_\lambda^{+}=\dfrac{S_\lambda+|S_\lambda|}{2}$ and Theorem 5 (ii).

(ii) Let $B^{\ast}=V(BB^{\ast})^{\frac{1}{2}}$ be the polar decomposition of $B^*,$ where $V$ is a partial isometry.
Setting $$\widetilde{Q}:=\frac{1}{2}\left(\begin{array}{cc}
\frac{T+T^{-1}+2I}{2}&(I+T^{-1})|B^{\ast}|\\ |B^{\ast}|(I+T^{-1})&2BB^{\ast}T^{-1}
\end{array}\right):
\mathcal{H}\oplus\mathcal{H}\rightarrow \mathcal{H}\oplus\mathcal{H},$$ and
$$\widetilde{X}:=\left(\begin{array}{cc}
(I+T)^{2}&2(I+T)|B^{\ast}|
\\2(I+T)|B^{\ast}|&4BB^{\ast}\end{array}\right):
\mathcal{H}\oplus\mathcal{H}\rightarrow \mathcal{H}\oplus\mathcal{H}
,$$
we get that
\begin{equation}S_1^{+}=\left(\begin{array}{cc}I&0\\0&V\end{array}\right)
\widetilde{Q}
\left(\begin{array}{cc}I&0\\0&V^{\ast}\end{array}\right),\end{equation}
and \begin{equation}\left(\begin{array}{cc}I&0\\0&V^{\ast}\end{array}\right)
\left(\begin{array}{cc}I&0\\0&V\end{array}\right)\widetilde{Q}=\widetilde{Q}.\end{equation}
On the other hand, it is clear that $T^{-1}|B^{\ast}|=|B^{\ast}|T^{-1}$ and $$\widetilde{Q}=\left(\begin{array}{cc}I+T^{2}+2T&2(I+T)|B^{\ast}|\\2(I+T)|B^{\ast}|&4BB^{\ast}\end{array}\right)
\left(\begin{array}{cc}\frac{1}{4}T^{-1}&0\\0&\frac{1}{4}T^{-1}\end{array}\right)=\widetilde{X}
\left(\begin{array}{cc}\frac{1}{4}T^{-1}&0\\0&\frac{1}{4}T^{-1}\end{array}\right),$$
which implies $R(\widetilde{Q})=R(\widetilde{X}).$
Thus Lemma 8 yields that
$$P_{\widetilde{Q}}=P_{\widetilde{X}}=\left(\begin{array}{cc}\frac{I+T^{-1}}{2}&T^{-1}|B^{\ast}|
\\|B^{\ast}|T^{-1}&2BB^{\ast}T^{-1}(I+T)^{-1}\end{array}\right).$$
Combining Lemma 9, equations (3.3) and (3.4), we have $$
\begin{array}{rl}P_{S_1^+}&=
\left(\begin{array}{cc}I&0\\0&V\end{array}\right)P_{\widetilde{Q}}
\left(\begin{array}{cc}I&0\\0&V^{\ast}\end{array}\right)\\\\&=\left(\begin{array}{cc}
\frac{I+T^{-1}}{2}&T^{-1}B\\B^{\ast}T^{-1}&2B^{\ast}T^{-1}(I+T)^{-1}B\end{array}\right). \end{array}$$

(iii) If $\lambda>0,$ then it is clear that $$S_{\lambda}^{+}
=\lambda\left(\begin{array}{cc}I&\frac{B}{\lambda}\\ \frac{B^{\ast}}{\lambda}&0\end{array}\right)^+=\lambda\left(\begin{array}{cc}I&B_0\\ B_0^{\ast}&0\end{array}\right)^+,$$
where $B_0=\frac{B}{\lambda}.$
 Thus (ii) implies that
$$T=(I+4B_0B_0^{\ast})^{\frac{1}{2}}=(I+\frac{4BB^{\ast}}{\lambda^2})^{\frac{1}{2}}
=\frac{(\lambda^2 I+4BB^{\ast})^{\frac{1}{2}}}{\lambda}$$ and $$P_{S_\lambda^+}=\left(\begin{array}{cc}\frac{I+T^{-1}}{2}&T^{-1}B_0\\B_0^{\ast}T^{-1}
 &2B_0^{\ast}T^{-1}(I+T)^{-1}B_0\end{array}
\right)=\left(\begin{array}{cc}
\frac{I+\lambda T_\lambda^{-1}}{2}
&T_\lambda^{-1}B\\ B^{\ast}T_\lambda^{-1}&
2B^{\ast}T_\lambda^{-1}(\lambda I+T_\lambda)^{-1}B\end{array}
\right).$$

(iv) If $\lambda<0,$ then $$S_\lambda^{+}=\dfrac{S_\lambda+|S_\lambda|}{2}
=\dfrac{\lambda}{2}[\left(\begin{array}{cc}I&\frac{B}{\lambda}\\ \frac{B^{\ast}}{\lambda}&0\end{array}\right)-\begin{vmatrix}\left(\begin{array}{cc}I&\frac{B}{\lambda}\\ \frac{B^{\ast}}{\lambda}&0\end{array}\right)\end{vmatrix}]=-\lambda\left(\begin{array}{cc}I&\frac{B}{\lambda}\\ \frac{B^{\ast}}{\lambda}&0\end{array}\right)^-.$$
Let $B_0=\frac{B}{\lambda}$ and $V_0$ is the unique partial isometry such that
$B_0^*=V_0(B_0B_0^{\ast})^{\frac{1}{2}}$ with $R(V_0)=\overline{R(B_0^*)}=\overline{R(B^*)}$ and $R(V_0^*)=\overline{R(B_0)}=\overline{R(B)}.$ Denoting
 $$S':=\left(\begin{array}{cc}I&\frac{B}{\lambda}\\ \frac{B^{\ast}}{\lambda}&0\end{array}\right)^-=\left(\begin{array}{cc}I&B_0\\ B_0^*&0\end{array}\right)^-,$$
 we conclude from the following Proposition 11 that
$$P_{S'}=\left(\begin{array}{cc}\frac{I-T^{-1}}{2}&-T^{-1}B_0\\-B_0^{\ast}T^{-1}
 &\frac{V_0(T^{-1}+I)V_0^*}{2}\end{array}
\right):\mathcal{H}\oplus\mathcal{K}\rightarrow \mathcal{H}\oplus\mathcal{K},$$  where $T=(I+4B_0B_0^{\ast})^{\frac{1}{2}}=
\frac{(\lambda^2 I+4BB^{\ast})^{\frac{1}{2}}}{-\lambda}.$
Clearly, $V=-V_0,$ then $$P_{S_\lambda^{+}}=P_{S'}=\left(\begin{array}{cc}\frac{I-T^{-1}}{2}&-T^{-1}B_0\\-B_0^{\ast}T^{-1}
 &\frac{V_0(T^{-1}+I)V_0^*}{2}\end{array}
\right)=\left(\begin{array}{cc}
\frac{I+\lambda T_\lambda^{-1}}{2}
&T_\lambda^{-1}B\\ B^{\ast}T_\lambda^{-1}
&\frac{V(I-\lambda T_\lambda^{-1})V^*}{2}\end{array}
\right).$$

(v) It is clear that with respect to the decomposition $\mathcal{H}\oplus\mathcal{H}:$ $\widetilde{C}:=\dfrac{1}{2}\left(\begin{array}{cc}
|B^*|&|B^*|\\ |B^*|&|B^*|\end{array}\right)\geqslant 0.$
By a direct calculation, we get that $$S_0^+=\frac{1}{2}\left(\begin{array}{cc} |B^*|&B\\ B^*&|B|\end{array}
\right)=\left(\begin{array}{cc}I&0\\0&V\end{array}\right)
 \widetilde{C}
\left(\begin{array}{cc}I&0\\0&V^{\ast}\end{array}\right)\ \ \hbox{ and }\ \ \left(\begin{array}{cc}I&0\\0&V^{\ast}V\end{array}\right)
\widetilde{C}
= \widetilde{C}.$$ Thus Lemma 9 implies $$\begin{array}{rcl}P_{S_0^+}&=&\left(\begin{array}{cc}I&0\\0&V\end{array}\right)
P_{\widetilde{C}}\left(\begin{array}{cc}I&0\\0&V^*\end{array}\right)
\\&=&\dfrac{1}{2}\left(\begin{array}{cc}I&0\\0&V\end{array}\right)
\left(\begin{array}{cc}P_{B}
&P_{B}\\ P_{B}&P_{B}\end{array}\right) \left(\begin{array}{cc}I&0\\0&V^*\end{array}\right)
\\&=&\dfrac{1}{2}\left(\begin{array}{cc}P_B&V^*\\
V&P_{B^*}\end{array}
\right).\end{array}$$

\ \ \ \ \ $\square$

{\bf Proposition 11.} Let $B\in \mathcal{B(K,H)}$ and $S\in \mathcal{B}(\mathcal{H}\oplus\mathcal{K})$ has the operator matrix form as $S_1$ in (3.1).  Then

(i) $$S^{-}= \left(\begin{array}{cc}\frac{T+T^{-1}-2I}{4}&\frac{(T^{-1}-I)B}{2}\\ \frac{B^{\ast}(T^{-1}-I)}{2}&B^{\ast}T^{-1}B\end{array}\right):
\mathcal{H}\oplus\mathcal{K}\rightarrow \mathcal{H}\oplus\mathcal{K},$$
where $T=(I+4BB^{\ast})^{\frac{1}{2}}.$

(ii) $P_{S}=diag(I, P_{B^*}):
\mathcal{H}\oplus\mathcal{K}\rightarrow \mathcal{H}\oplus\mathcal{K}.$

(iii)
$$P_{S^-}=\left(\begin{array}{cc}\frac{I-T^{-1}}{2}&-T^{-1}B
\\-B^{\ast}T^{-1}& \frac{V(T^{-1}+I)V^*}{2}\end{array}\right):
\mathcal{H}\oplus\mathcal{K}\rightarrow \mathcal{H}\oplus\mathcal{K},$$ where $V$ is the unique partial isometry such that
$B^{\ast}=V(BB^{\ast})^{\frac{1}{2}}$ with $R(V)=\overline{R(B^*)}$ and $R(V^*)=\overline{R(B)}.$

{\bf Proof.} (i) is immediate from the fact $S^{-}=S^{+}-S.$

(ii) Let $V$ be the unique partial isometry such that
$B^{\ast}=V(BB^{\ast})^{\frac{1}{2}},$ $R(V)=\overline{R(B^*)}$ and $R(V^*)=\overline{R(B)}.$ Setting $$\widetilde{V}=diag(I,V):\mathcal{H}\oplus\mathcal{K}\rightarrow \mathcal{H}\oplus\mathcal{K}$$
and
 $$Q=\left(\begin{array}{cc}I&|B^*|\\|B^{\ast}|& 0\end{array}\right):\mathcal{H}\oplus\mathcal{H}\rightarrow \mathcal{H}\oplus\mathcal{H},$$ we claim that
$N(Q)=\{0\}\oplus N(B^*).$ Indeed, suppose that
$x\in \mathcal{H}$ and $y\in \mathcal{H}$ satisfy $Q(x\oplus y)=0.$ Then
$$x+|B^*|y=0 \ \ \ \ \hbox{ and } \ \ \ \ |B^*|x=0,$$ so $BB^*y=|B^*|(x+|B^*|y)=0,$ which implies
$y\in N(B^*)$ and $x=0.$
Thus $$P_{Q}=diag(I, P_{B}):\mathcal{H}\oplus\mathcal{H}\rightarrow \mathcal{H}\oplus\mathcal{H}.$$
Clearly, $S=\widetilde{V}Q\widetilde{V}^*$ and $\widetilde{V}^*\widetilde{V}Q=Q,$
so Lemma 9 yields that  $$P_{S}=\widetilde{V}P_{Q}\widetilde{V}^*=diag(I, VP_{B}V^*)=diag(I, VV^*)=diag(I, P_{B^*}).$$

(iii) It is easy to verify that
$$\begin{array}{rcl}I-2BB^*T^{-1}(I+T)^{-1}&=&(T+T^2)^{-1}(T+T^2-2BB^*)
\\&=&\dfrac{(T+T^2)^{-1}(2I+4BB^*+2T)}{2}\\&=&\dfrac{(T+T^2)^{-1}(I+T)^2}{2}\\&=&\dfrac{T^{-1}+I}{2}.\end{array}$$ Hence
$$P_{B^*}-2B^{\ast}T^{-1}(I+T)^{-1}B=VV^*-2VBB^*T^{-1}(I+T)^{-1}V^*=\frac{V(T^{-1}+I)V^*}{2},$$ which implies
$$P_{S^-}=P_{S}-P_{S^+}=\left(\begin{array}{cc}\frac{I-T^{-1}}{2}&-T^{-1}B
\\-B^{\ast}T^{-1}& \frac{V(T^{-1}+I)V^*}{2}\end{array}\right):
\mathcal{H}\oplus\mathcal{K}\rightarrow \mathcal{H}\oplus\mathcal{K}.$$\ \ \ \ \ $ \square$

The following is a direct corollary from Theorem 10.

{\bf Corollary 12.} Let $E\in {\mathcal{B(H)}}^{Id}$ have the form (1.1). Then with respect to the decomposition $\mathcal{H}=R(E)\oplus R(E)^{\perp},$

(i)$$(E+E^{\ast})^{+}=\dfrac{1}{2}\left(\begin{array}{cc}T+T^{-1}+2I&(I+T^{-1})E_1
\\E_1^{\ast}(I+T^{-1})&E_1^{\ast}T^{-1}E_1\end{array}\right),$$

(ii)$$P_{(E+E^{\ast})^+}=\dfrac{1}{2}\left(\begin{array}{cc}I+T^{-1}&T^{-1}E_1
\\E_1^{\ast}T^{-1}&E_1^{\ast}(T+T^{2})^{-1}E_1\end{array}\right),$$
where $T=(I+E_1E_1^{\ast})^{\frac{1}{2}}.$

\section{$J$-projections}

A projection $E\in \mathcal{B(H)}^{Id}$ is said to be $J$-projection if $E=JE^{\ast}J,$
which is equivalent to say that $JE$ is self-adjoint.
A $J$-projection $E$ is called to be positive (negative) if $JE\geqslant0 (JE\leqslant0$).
In this section,
 we mainly characterize the symmetry $J$ such that a projection $E$ is the $J$-projection.
In particular, the minimal element of the set of all symmetries $J$ with $JE\geqslant 0$ is given.
The following lemma is needed.

{\bf Lemma 13.} Let $E\in \mathcal{B(H)}^{Id}$ have the form (1.1). If $J$ is a symmetry,
 then $E$ is a $J$-projection if and only if $J$ has the operator matrix form
\begin{equation}J=\left(\begin{array}{cc}J_1(I+E_1E_1^{\ast})^{-\frac{1}{2}}&J_1(I+E_1E_1^{\ast})^{-\frac{1}{2}}E_1\\
E_1^{\ast}(I+E_1E_1^{\ast})^{-\frac{1}{2}}
J_1&J_2(I+E_1^{\ast}E_1)^{-\frac{1}{2}}
\end{array}\right):R(E)\oplus R(E)^{\perp}\rightarrow R(E)\oplus R(E)^{\perp},\end{equation}
 where $J_1$ and $J_2$ are  symmetries on the subspaces
 $R(E)$ and $R(E)^{\perp}$ respectively, satisfying $J_1E_1+E_1J_2=0.$

{\bf Proof.} Sufficiency. If $J_1E_1+E_1J_2=0,$ then $J_1E_1=-E_1J_2,$ which induces
\begin{equation}J_1E_1E_1^{\ast}=-E_1J_2E_1^{\ast}=-E_1(-E_1^{\ast}J_1)=E_1E_1^{\ast}J_1.\end{equation}
Similarly, we get that \begin{equation}J_2E_1^{\ast}E_1=E_1^{\ast}E_1J_2.\end{equation}
 Combining equations (1.1) and (4.1)-(4.3), we easily verify the equation $JE^{\ast}J=E.$

Necessity. Without loss of generality, we may assume
\begin{equation}J=\left(\begin{array}{cc}
Q_{11}&Q_{12}\\Q_{12}^{\ast}&Q_{22}
\end{array}\right):R(E)\oplus R(E)^{\perp}\rightarrow R(E)\oplus R(E)^{\perp},\end{equation}
where $Q_{11}$ and $Q_{22}$ are self-adjoint operators.
It follows from the fact $J^2=I$  that \begin{equation}\begin{cases}Q_{11}^2+
Q_{12}Q_{12}^{\ast}=I\ \qquad\qquad\textcircled{1}
&\\Q_{11}Q_{12}+Q_{12}Q_{22}=0\ \ \qquad \textcircled{2}&\\Q_{12}^{\ast}Q_{12}+Q_{22}^2=I. \qquad\qquad\textcircled{3}\end{cases}\end{equation}

On the other hand, $E=JE^{\ast}J$ implies that $JE$ is self-adjoint, so
\begin{equation}Q_{12}=Q_{11}E_1.\end{equation}
Using equations $ \textcircled{1}$ of (4.5) and (4.6), we have $$Q_{11}(I+E_1E_1^{\ast})Q_{11}=I,$$
which yields $Q_{11}$ is invertible on the subspace $R(E)$ and $Q_{11}^2=(I+E_1E_1^{\ast})^{-1}.$
Thus $$|Q_{11}|=(I+E_1E_1^{\ast})^{-\frac{1}{2}}$$ follows from the fact that  $Q_{11}$ is self-adjoint.
Let $Q_{11}=V|Q_{11}|$ be the polar decomposition of $Q_{11}.$
 Then $V=Q_{11}|Q_{11}|^{-1}$ satisfies
$V=V^{\ast}=V^{-1},$  so $J_1:=V$ is a symmetry on the subspace $R(E)$ and \begin{equation}Q_{11}=J_1|Q_{11}|=J_1(I+E_1E_1^{\ast})^{-\frac{1}{2}}.\end{equation}
Similarly,  equations $ \textcircled{3}$ of (4.5) and (4.6) imply
$$Q_{22}^2=I-E_1^*Q_{11}^2E_1=I-E_1^*(I+E_1E_1^{\ast})^{-1}E_1=
I-(I+E_1^*E_1)^{-1}E_1^{\ast}E_1=(I+E_1^{\ast}E_1)^{-1},$$ since
$(I+E_1^*E_1)^{-1}E_1^{\ast}=E_1^{\ast}(I+E_1E_1^*)^{-1}.$  Thus
\begin{equation}|Q_{22}|=(I+E_1^{\ast}E_1)^{-\frac{1}{2}}.\end{equation}
Setting $J_2:=Q_{22}|Q_{22}|^{-1},$  we easily verify  that $J_2=J_2^{\ast}=J_2^{-1}$  and
\begin{equation}Q_{22}=J_2|Q_{22}|=J_2(I+E_1^{\ast}E_1)^{-\frac{1}{2}}.\end{equation}
Moreover, by equations $ \textcircled{2}$ of (4.5) and (4.7)-(4.9), we get that $$(J_1E_1+E_1J_2)(I+E_1^{\ast}E_1)^{-\frac{1}{2}}=Q_{11}E_1+E_1Q_{22}=0,$$
which induces $J_1E_1+E_1J_2=0.$  \qquad $\square$

In the following, we give some new characterizations for the positivity of the $J$-projections.

{\bf Corollary 14.} Let $E\in \mathcal{B(H)}^{Id}$ have the form (1.1). Then $JE\geqslant0$
if and only if $$J=\left(\begin{array}{cc}(I+E_1E_1^{\ast})^{-\frac{1}{2}}&(I+E_1E_1^{\ast})^{-\frac{1}{2}}E_1\\
E_1^{\ast}(I+E_1E_1^{\ast})^{-\frac{1}{2}}
&J_2(I+E_1^{\ast}E_1)^{-\frac{1}{2}}
\end{array}\right):R(E)\oplus R(E)^{\perp}\rightarrow R(E)\oplus R(E)^{\perp},$$
 where $J_2$ is a symmetry on the subspace
 $R(E)^{\perp}$  with $E_1=-E_1J_2.$

{\bf Proof.} Sufficiency. With respect to the decomposition $R(E)\oplus R(E)^{\perp},$ it is clear that $$JE=\left(\begin{array}{cc}(I+E_1E_1^{\ast})^{-\frac{1}{2}}&(I+E_1E_1^{\ast})^{-\frac{1}{2}}E_1\\
E_1^{\ast}(I+E_1E_1^{\ast})^{-\frac{1}{2}}&
E_1^{\ast}(I+E_1E_1^{\ast})^{-\frac{1}{2}}
E_1\end{array}\right).$$
Calculating the Schur complement of $JE,$ we know $JE\geqslant0.$

 Necessity. By Lemma 13, we only need to show $J_1=I.$
 It is easy to see that $JE\geqslant0$ implies
 $$Q_{11}=J_1|Q_{11}|=J_1(I+E_1E_1^{\ast})^{-\frac{1}{2}}\geqslant0,$$ so
  $$J_1(I+E_1E_1^{\ast})^{-\frac{1}{2}}=Q_{11}=|Q_{11}|=(I+E_1E_1^{\ast})^{-\frac{1}{2}},$$
which yields $J_1=I$ as desired. \qquad $\square$

{\bf Theorem 15.} Let $E\in \mathcal{B(H)}^{Id}.$ Then  \begin{equation}min\{J: JE\geqslant0,\ J=J^{\ast}=J^{-1}\}=2P_{(E+E^{\ast})^{+}}-I,\end{equation}
where the ``min" is in the sense of Loewner partial order.

{\bf Proof.} Suppose that $E\in \mathcal{B(H)}^{Id}$ has the form (1.1). Then by Corollary 12,
 we get that $$2P_{(E+E^{\ast})^{+}}-I=
\left(\begin{array}{cc}T^{-1}&T^{-1}E_1
\\E_1^{\ast}T^{-1}&E_1^{\ast}(T+T^2)^{-1}
E_1-I\end{array}\right):R(E)\oplus R(E)^{\perp}\rightarrow R(E)\oplus R(E)^{\perp},$$
where $T=(I+E_1E_1^{\ast})^{\frac{1}{2}}.$

It is easy to see that  $$\begin{array}{rl}E_1^{\ast}(T+T^2)^{-1}E_1-I&=E_1^{\ast}[(I+E_1E_1^{\ast})^{\frac{1}{2}}+I+E_1E_1^{\ast}]^{-1}E_1-I
\\&=[(I+E_1^{\ast}E_1)^{\frac{1}{2}}+I+E_1^{\ast}E_1]^{-1}E_1^{\ast}E_1-I
\\&=-[(I+E_1^{\ast}E_1)^{\frac{1}{2}}+I+E_1^{\ast}E_1]^{-1}[I+(I+E_1^{\ast}E_1)^{\frac{1}{2}}]
\\&=-(I+E_1^{\ast}E_1)^{-\frac{1}{2}}.\end{array}$$

Setting $J_2=-I$ in proof of Corollary 14, we get $(2P_{(E+E^{\ast})^{+}}-I)E\geqslant0,$
so the left hand is no more than the right hand of equation (4.10).

 If $JE\geqslant0,$ then by Corollary 14, we have $$J=\left(\begin{array}{cc}(I+E_1E_1^{\ast})^{-\frac{1}{2}}&(I+E_1E_1^{\ast})^{-\frac{1}{2}}E_1\\
E_1^{\ast}(I+E_1E_1^{\ast})^{-\frac{1}{2}}&
J_2(I+E_1^{\ast}E_1)^{-\frac{1}{2}}\end{array}
\right):R(E)\oplus R(E)^{\perp}\rightarrow R(E)\oplus R(E)^{\perp}.$$  Obviously, $J_2E_1^*E_1=E_1^*E_1J_2$ implies $$J_2(I+E_1^{\ast}E_1)^{-\frac{1}{2}}-[-(I+E_1^{\ast}E_1)^{-\frac{1}{2}}]=(I+J_2)(I+E_1^{\ast}E_1)^{-\frac{1}{2}}\geqslant0,$$
which yields $J\geqslant2P_{(E+E^{\ast})^{+}}-I$ as desired.\qquad $\square$

{\bf Proposition 16.}([11, Proposition 4]) Let $E\in \mathcal{B(H)}^{Id}$ and $J$ be a symmetry.  If $E=JE^{\ast}J,$ then

(a) $JE\geqslant0$ if and only if $J\geqslant2P_{(E+E^{\ast})^{+}}-I.$

(b) $JE\geqslant0$ if and only if $(J+I)E(J+I)\geqslant0  \ \ \ \hbox{ and } \ \ \ R(E)\cap R(I-J)=\{0\}.$

{\bf Proof.} (a) Necessity is clear by Theorem 15.

Sufficiency.  Since $E=JE^{\ast}J,$  it follows from Lemma 13 that $$J=\left(\begin{array}{cc}J_1(I+E_1E_1^{\ast})^{-\frac{1}{2}}&J_1(I+E_1E_1^{\ast})^{-\frac{1}{2}}E_1\\
E_1^{\ast}(I+E_1E_1^{\ast})^{-\frac{1}{2}}J_1
&J_2(I+E_1^{\ast}E_1)^{-\frac{1}{2}}
\end{array}\right):R(E)\oplus R(E)^{\perp}\rightarrow R(E)\oplus R(E)^{\perp}.$$ Then $J\geqslant2P_{(E+E^{\ast})^{+}}-I$ yields
$J_1(I+E_1E_1^{\ast})^{-\frac{1}{2}}\geqslant(I+E_1E_1^{\ast})^{-\frac{1}{2}},$ so
$J_1\geqslant I$ follows from the fact $J_1E_1E_1^*=E_1E_1^*J_1.$
Thus $J_1=I,$ as $J_1=J_1^*=J_1^{-1}.$  Then Corollary 14 induces $JE\geqslant0$ as desired.

 (b)  Necessity. If $JE\geqslant0,$ then $$(J+I)E(J+I)=(J+I)JE(J+I)\geqslant0.$$
Let $x\in R(E)\cap R(I-J).$ Then there exists $y\in \mathcal{H}$ such that $x=Ex=(I-J)y,$  so $$Jx=JEx=J(I-J)y=(J-I)y=-x,$$
which yields $$0\leqslant\langle JEx,x\rangle=\langle -x,x\rangle=-\|x\|^{2}.$$
 Thus $x=0.$

Sufficiency. By Corollary 14, we need to prove that $J_1=I.$

Since $J_1=J_1^*=J_1^{-1},$ then without loss of generality,
we assume \begin{equation}J_1=\left(\begin{array}{cc}
I&0\\0&-I\end{array}\right):\mathcal{M}\oplus \mathcal{M}^{\perp}\rightarrow \mathcal{M}\oplus \mathcal{M}^{\perp},\end{equation}
 where $\mathcal{M}\subseteq R(E)$ and $\mathcal{M}^{\perp}=R(E)\ominus \mathcal{M}.$
Then $J_1E_1E_1^{\ast}=E_1E_1^{\ast}J_1$ implies that $\mathcal{M}$ is a reduced subspace of $E_1E_1^{\ast},$
so \begin{equation}(I+E_1E_1^{\ast})^{-\frac{1}{2}}
=\left(\begin{array}{cc}Q_1&0\\0&Q_2
\end{array}\right):\mathcal{M}\oplus \mathcal{M}^{\perp}\rightarrow \mathcal{M}\oplus \mathcal{M}^{\perp},\end{equation}
 where $Q_1$ and $Q_2$ are positive and invertible.
Using equation (4.7), we have \begin{equation}Q_{11}=J_1(I+E_1E_1^{\ast})^{-\frac{1}{2}}
=\left(\begin{array}{cc}Q_1&0\\0&-Q_2
\end{array}\right):\mathcal{M}\oplus \mathcal{M}^{\perp}\rightarrow \mathcal{M}\oplus \mathcal{M}^{\perp}.\end{equation}
  By equations (1.1),(4.1) and a direct calculation, we know that $(J+I)E(J+I)\geqslant0$ implies $$(Q_{11}+I)^2+(Q_{11}+I)E_1E_1^{\ast}Q_{11}\geqslant0,$$
so $$Q_{11}(I+E_1E_1^*)Q_{11}+Q_{11}+(I+E_1E_1^*)Q_{11}+I=(Q_{11}+I)^2+(Q_{11}+I)E_1E_1^{\ast}Q_{11}\geqslant0.$$ Using again equation (4.7), we get that  $$Q_{11}(I+E_1E_1^*)Q_{11}=I  \hbox{   }\hbox{ and }\hbox{   }
(I+E_1E_1^*)Q_{11}=Q_{11}^{-1},$$ so $$Q_{11}+2I+Q_{11}^{-1}
=Q_{11}(I+E_1E_1^*)Q_{11}+Q_{11}+
(I+E_1E_1^*)Q_{11}+I\geqslant0.$$
Then
$$\left(\begin{array}{cc}Q_1+2I+Q_1^{-1}&0\\0&-Q_2+2I-Q_2^{-1}\end{array}\right)\geqslant0$$ follows from equation (4.13),
which yields $$-(Q_2^\frac{1}{2}-Q_2^{-\frac{1}{2}})^2=-Q_2+2I-Q_2^{-1}\geqslant0,$$
so $Q_2^\frac{1}{2}-Q_2^{-\frac{1}{2}}=0.$ Thus $Q_2=I,$ which implies
Thus \begin{equation} (I+E_1E_1^{\ast})^{-\frac{1}{2}}=\left(\begin{array}{cc}Q_1&0\\0&I\end{array}\right) \ \ \
\hbox{ and } \ \ \ \ Q_{11}=\left(\begin{array}{cc}Q_1&0\\0&-I\end{array}\right).\end{equation}

Let $x\in \mathcal{M}^{\perp}.$ Then $(I+E_1E_1^{\ast})^{-\frac{1}{2}}x=x,$
which yields $(I+E_1E_1^{\ast})x=x,$ so $E_1E_1^{\ast}x=0.$
Hence, \begin{equation} x\in N(E_1^{\ast})\subseteq R(E).\end{equation}

Moreover, combining equations (4.1) and (4.14), we conclude that
$$Jx=\left(\begin{array}{cc}Q_{11}x\\E_1^{\ast}Q_{11}x\end{array}\right)
=\left(\begin{array}{cc}-x\\-E_1^{\ast}x\end{array}\right)
=\left(\begin{array}{cc}-x\\0\end{array}\right),$$ which induces $(I-J)x=x-Jx=2x,$
so \begin{equation}x=(I-J)\dfrac{x}{2}\in R(I-J).\end{equation}
Thus equations (4.15) and (4.16) imply $$x\in R(I-J)\cap R(E)=\{0\}.$$
 Then $\mathcal{M}^{\perp}=0,$ so $ J_1=I$ as desired.
 \ \ \ \ \ $\square$

{\bf Proposition 17.} ([2, Theorem 2.1]) Let $\mathcal{M}$ be a non-trivial closed subspace of $\mathcal{H}$ and $J$ be a symmetry.
Then there exists $E\in\mathcal{B(H)}^{Id}$ such that $R(E)=\mathcal{M}$ and
$JE^{\ast}J=E$ if and only if $P_{\mathcal{M}}J\mid_{\mathcal{M}}$
is invertible on the subspace $\mathcal{M}.$
In this case, $E$ is unique.

{\bf Proof.} Let $$J=\left(\begin{array}{cc}Q_{11}&Q_{12}
\\Q_{12}^{\ast}&Q_{22}\end{array}\right):
\mathcal{M}\oplus \mathcal{M}^{\perp}\rightarrow \mathcal{M}\oplus \mathcal{M}^{\perp},$$
where $Q_{11}$ and $Q_{22}$ are self-adjoint operators.

Sufficiency. If $P_{\mathcal{M}}J\mid_{\mathcal{M}}$ is invertible on the subspace $\mathcal{M},$
then $Q_{11}$ is invertible.
Define an idempotent operator $E$ as $$E=\left(\begin{array}{cc}I& Q_{11}^{-1}Q_{12}\\0&0\end{array}\right):
\mathcal{M}\oplus \mathcal{M}^{\perp}\rightarrow\mathcal{M}\oplus \mathcal{M}^{\perp}.$$
It is easy to verify that $R(E)=\mathcal{M}$ and
$JE^{\ast}J=E.$

Necessity. Suppose that there exists $E\in\mathcal{B(H)}^{Id}$ such that $R(E)=\mathcal{M}$ and
$JE^{\ast}J=E.$ Then by Lemma 13 $$J=\left(\begin{array}{cc}
J_1(I+E_1E_1^{\ast})^{-\frac{1}{2}}
&J_1(I+E_1E_1^{\ast})^{-\frac{1}{2}}E_1\\
E_1^{\ast}(I+E_1E_1^{\ast})^{-\frac{1}{2}}
J_1&J_2(I+E_1^{\ast}E_1)^{-\frac{1}{2}}
\end{array}\right):\mathcal{M}\oplus \mathcal{M}^{\perp}\rightarrow \mathcal{M}\oplus \mathcal{M}^{\perp},$$ where $E$ has the form as (1.1). Thus $P_{\mathcal{M}}J\mid_{\mathcal{M}}=
J_1(I+E_1E_1^{\ast})^{-\frac{1}{2}}$
is invertible on the subspace $\mathcal{M}.$
Also, Lemma 13 (the above equation) implies $E_1=(P_{\mathcal{M}}J\mid_{\mathcal{M}})^{-1}
P_{\mathcal{M}}J\mid_{\mathcal{M}^{\perp}}$
is unique, which says the uniqueness of $E.$    \ \ \ \ \ $\square$

The following result is obtained in [2,7,11].
However, our methods are completely different from those of [2,7,11].

{\bf Corollary 18.} For a $J$-projection $E$, there exists uniquely a $J$-positive projection $Q,$
a $J$-negative projection $R$ such that $$E=Q+R,\ \ \  QR=RQ=0\ \ \ \hbox {and} \ \ \ QR^{\ast}=R^{\ast}Q=0.$$

{\bf Proof.} Suppose that $E$ has operator matrix form (1.1).
Using Lemma 13, we conclude that $JEJ=E^{\ast}$ is equivalent to that $J$ has the form (4.1).
With respect to the space decomposition $\mathcal{H}=R(E)\oplus R(E)^{\perp},$ we define $Q$ and $R$ as \begin{equation}Q:=\left(\begin{array}{cc}
\dfrac{I_1+J_1}{2}&\dfrac{I_1+J_1}{2}E_1
\\0&0\end{array}\right) \quad \hbox{ and }\  \  \
R:=\left(\begin{array}{cc}\dfrac{I_1-J_1}{2}&
\dfrac{I_1-J_1}{2}E_1\\\\0&0\end{array}\right),\end{equation}
where $I_1$ is the identity operator
on the subspace $R(E).$  Clearly, $Q$ and $R$ are projections with $E=Q+R.$
A direct calculation yields that $QR=RQ=0\ \hbox{and}\ QR^{\ast}=R^{\ast}Q=0.$
Moreover, combining equations (4.1) and (4.17), we get that $$JQ=\left(\begin{array}{cc}I_1&0\\0&E_1^{\ast}\end{array}\right)
\left(\begin{array}{cc}\dfrac{I_1+J_1}{2}T^{-1}&
\dfrac{I_1+J_1}{2}T^{-1}\\\\
\dfrac{I_1+J_1}{2}T^{-1}&\dfrac{I_1+J_1}{2}T^{-1}\end{array}\right)
\left(\begin{array}{cc}I_1&0\\0&E_1\end{array}\right)\geqslant0$$
and $$JR=-\left(\begin{array}{cc}I_1&0\\0&E_1^{\ast}\end{array}\right)
\left(\begin{array}{cc}\dfrac{I_1-J_1}{2}T^{-1}
&\dfrac{I_1-J_1}{2}T^{-1}\\\\
\dfrac{I_1-J_1}{2}T^{-1}&\dfrac{I_1-J_1}{2}T^{-1}\end{array}\right)
\left(\begin{array}{cc}I_1&0\\0&E_1\end{array}\right)\leqslant0,$$ since $\dfrac{I_1+J_1}{2}T^{-1}\geqslant0$ and $\dfrac{I_1-J_1}{2}T^{-1}\geqslant0,$ where $T=(I+E_1E_1^*)^{\frac{1}{2}}.$
Thus $Q$ is a $J$-positive projection and
 $R$ is a $J$-negative projection.

 To show uniqueness, we assume that $E=Q_1+Q_2$, where $Q_1$ and $Q_2$ satisfy $JQ_1\geqslant0,$\ $JQ_2\leqslant0$
and $Q_1Q_2^{\ast}=Q_2^{\ast}Q_1=0.$ It is obvious that $$JE=JQ_1+JQ_2 \ \
\hbox{ and }\ \ (JQ_2) (JQ_1)=(Q_2^{\ast}J)(JQ_1)=Q_2^{\ast}Q_1=0.$$
 Since $(JR)(JQ)=0$ and $JE=JQ+JR,$ then $JQ=JQ_1$ and $ JR=JQ_2$ follow from the uniqueness of the positive part and the negative part for a self-adjoint operator.
 Thus $Q=Q_1$ and $ R=Q_2.$ \qquad  $\square$

{\bf Remark 19.} In general, a $J$-projection $E$ is representable (not uniquely) as $E=E_++E_-,$
where $E_+$ and $E_-$ are positive and negative $J$-projection, respectively.


\begin{thebibliography}{99}
\bibitem {s1} T. Ya. Azizov and I.S. Iokhvidov, Linear operators in spaces with an indefinite metric, John Wiley and Sons, 1989.
\bibitem {s1} T. Ando, Projections in Krein spaces, Linear Algebra Appl. 12 (2009), 2346-2358.
\bibitem {s1} M. L. Arias, G. Corach,  A. Maestripieri, Products of Idempotent Operators,
Integr. Equ. Oper. Theory 88 (2017), 269-286.
\bibitem {s1} R. Bhatia, Matrix Analysis, Graduate Texts in Mathematics, vol.169, Springer-Verlag, New York (1997).
\bibitem {s1} J. B. Conway, A course in operator theory, Graduate studies in mathematics, 21 (2000).
\bibitem {s1} G. Corach, A. Maestripieri and D. Stojanoff, Oblique projections and Schur
complements, Acta Sci. Math. (Szeged) 67 (2001), 337-356.
\bibitem {s1} S. Hassi, K. Nordstrom, On projections in a space with an indefinite metric, Linear Algebra Appl. 208/209 (1994) 401-407.
\bibitem {s1} A. Maestripieri, F. M. Per\'{\i}a,  Normal Projections in Krein Spaces, Integr. Equ. Oper. Theory 76 (2013), 357-380.
\bibitem {s1} A. Maestripieri, F. M. Per\'{\i}a,  Decomposition of selfadjoint Projections in Krein Spaces, Acta Sci. Math. (Szeged) 72 (2006), 611-638.
\bibitem {s1} M. Matvejchuk, Idempotents in a space with conjugation, Linear Algebra and its Applications 438 (2013) ,71-79.
\bibitem {s1} M. Matvejchuk, Idempotents as J-Projections, Int. J. Theor Phys., 50 (2011), 3852-3856.

\bibitem {s1} M. Matvejchuk, Idempotents and Krein spaces, Lobachevskii J. Math., 32 (2) (2011), 128-134.






\end{thebibliography}
\end{document}